\documentclass[11pt,reqno]{amsart}
\allowdisplaybreaks[4]
\usepackage{amssymb}
\usepackage{anysize}
\usepackage{hyperref}
\usepackage{extarrows}

\numberwithin{equation}{section}
\newtheorem{thm}{Theorem}[section]

\newtheorem{rem}[thm]{Remark}

\newcommand{\be}{\begin{equation}}
\newcommand{\ee}{\end{equation}}
\newcommand{\bea}{\begin{eqnarray*}}
\newcommand{\eea}{\end{eqnarray*}}

\makeatletter

\newcommand{\Rmnum}[1]{\expandafter\@slowromancap\romannumeral #1@}
\makeatother
\makeatletter 
\@addtoreset{equation}{section}
\makeatother  

\marginsize{30mm}{30mm}{35mm}{35mm}

\begin{document}

\title[]{On asymptotic expansion and CLT of linear eigenvalue statistics for sample covariance matrices when $N/M\rightarrow0$$*$}
\author{Zhigang Bao}
\address{Department of Mathematics, Zhejiang University, Hangzhou, Zhejiang 310027, P.R. China}
\thanks{$*$The work was supported partially by NSFC grant
11071213, ZJNSF grant R6090034  and   SRFDP grant 20100101110001.  }
\email{maomie2007@gmail.com}

\dedicatory
 {
 Department of Mathematics, Zhejiang University
 Hangzhou, 310027, P. R. China.
 }


\date{\today}


\maketitle


\begin{abstract}
We study the renormalized real sample covariance matrix $H=X^TX/\sqrt{MN}-\sqrt{M/N}$ with $N/M\rightarrow0$ as $N, M\rightarrow \infty$ in this paper. And we always assume $M=M(N)$. Here $X=[X_{jk}]_{M\times N}$ is an $M\times N$ real random matrix with i.i.d entries, and we assume $\mathbb{E}|X_{11}|^{5+\delta}<\infty$ with some small positive $\delta$. The Stieltjes transform $m_N(z)=N^{-1}Tr(H-z)^{-1}$ and the linear eigenvalue statistics of $H$ are considered. We mainly focus on the asymptotic expansion of $\mathbb{E}\{m_N(z)\}$ in this paper. Then for some fine test function, a central limit theorem for the linear eigenvalue statistics of $H$ is established.  We show that the variance of the limiting normal distribution coincides with the case of a real Wigner matrix with Gaussian entries.

\end{abstract}

\maketitle

\section{Introduction }
As an important branch of the Random Matrix Theory, the study towards sample covariance matrix traces back to the work of Hsu \cite{Hsu} and Wishart \cite{WJ}. The modern formulation of a large dimensional sample covariance matrix always indicates the matrix in the form of  $S_N=\frac1MX^TX$, where $X=[X_{jk}]_{M\times N}$ is an $M\times N$ real random matrix with mean zero and variance $\sigma^2$ i.i.d entries and $N/M\rightarrow y\in(0,\infty)$ as $N, M:=M(N)$ tend to infinity. So usually the sample size $M$ and parameter number $N$ are assumed to be with the same order. For a symmetric matrix with eigenvalues $\lambda_1,\cdots,\lambda_N$, we denote its empirical spectral distribution by $F_N(x)=\frac{1}{N}\sum_{i=1}^{N}\mathbf{1}_{\{\lambda_i\leq x\}}$. And the limit of $F_N(x)$ as $N\rightarrow\infty$ is often referred to as limiting spectral distribution. In the case of $S_N$, it is well known that Mar$\check{c}$enko and Pastur \cite{MP} firstly found: as $N$ tends to infinity $F_N(x)$ almost surely converges to the so called Mar$\check{c}$enko-Pastur (MP) law with the density function
$$
\rho_y(x)=\frac{1}{2\pi xy\sigma^2}\sqrt{(b-x)(x-a)}\mathbf{1}_{\{a\leq x\leq b\}},
$$
and there is a point mass $1-1/y$ at the origin if $y>1$, where $a=\sigma^2(1-\sqrt{y})^2$ and $b=\sigma^2(1+\sqrt{y})^2$.

While in modern statistics, the case of $N/M\rightarrow0$ as both $N$ and $M$ tend to infinity is also very common, see El Karoui \cite{EK} for example. In this paper, we will focus on the theoretical aspect of such a particular case of the sample covariance matrix. It is easy to see when $N/M\rightarrow y=0$, the interval $[a,b]$ will shrink to the point $\sigma^2$, so every eigenvalue of $S$ will tend to $\sigma^2$ then. If we centralize $S$ by subtracting $\sigma^2I_N$ and multiplying it with $\sqrt{M/N}$, the range of the eigenvalues will be enlarged to be order $1$ typically. In fact, under the assumption $N/M\rightarrow 0$ as $N\rightarrow\infty$, it was understood long time ago that the matrix
\begin{eqnarray}
H_N:=\sqrt{\frac MN}\left(\frac1MX^TX-\sigma^2I_N\right) \label{1.1}
\end{eqnarray}
behaves similar with a Wigner matrix of dimension $N$ on many spectral properties.

A real Wigner matrix can be defined as a real symmetric random matrix $W_N$ with mean zero, finite variance i.i.d diagonal entries and mean zero, variance $\omega^2/N$ i.i.d above-diagonal entries, and all entries on or above diagonal are independent. As a cornerstone in the Random Matrix Theory, the so called semicircle law as the limiting spectral distribution of $W_N$ was firstly raised by Wigner \cite{WE}, with the density given by
\begin{eqnarray}
\rho_{sc}(x)=\frac{1}{2\pi\omega^2}\sqrt{4\omega^2-x^2}\mathbf{1}_{\{|x|\leq2\omega\}}.\label{1.2}
\end{eqnarray}

It is not difficult to see, when we fix $N$ and only let $M$ tends to infinity, $H_N$ will tend to a Wigner matrix with Gaussian entries under the effect of the central limit theorem.
Thus if we let both $N$ and $M$ tend to infinity with $N/M\rightarrow 0$, it is natural to ask whether $F_N(x)$ of $H_N$ tend to the semicircle law as well. A rigorous proof was given by Bai and Yin \cite{BY} through a moment method: under the assumption of the existence of $4$-th moment of $X_{jk}$, $F_N(x)$ converges almost surely to the semicircle law (\ref{1.2}) with $\omega=\sigma^2$. So a question naturally arises, how close are the spectral properties of the Wigner matrix $W_N$ and the renormalized sample covariance matrix $H_N$ when we set $\omega=\sigma^2$? However, such a question makes no sense if we do not specify some particular spectrum statistics to compare. So in this paper, we will study some representative objects for $H_N$, like Stieltjes transform of $F_N(x)$ and linear eigenvalue statistics of $H_N$. Then compare them to the corresponding results known for the Wigner matrix $W_N$.

The Stieltjes transform of $F_{N}(x)$ is defined by
\begin{eqnarray*}
m_N(z)=\int\frac{1}{x-z}dF_N(x) \label{1.3}
\end{eqnarray*}
for any $z=E+i\eta$ with $E\in\mathbb{R}$ and $\eta>0$. We just consider fixed $z$ , and specify $F_N(x)$ as the empirical spectral distribution of $H_N$ throughout the paper. If we introduce the Green function $G_N(z)=(H_N-z)^{-1}$, we also have
\begin{eqnarray*}
m_N(z)=\frac1N TrG_N(z)=\frac1N\sum_{j=1}^{N}G_{jj}. \label{1.4}
\end{eqnarray*}
Here we denote $G_{jk}$ as the $(j,k)$ entry of $G_N(z).$ Both as a powerful tool and as a relevant spectral statistic, the Stieltjes transform of the empirical spectral distribution has shown to be particularly important in the Random Matrix Theory. As is well known, the convergence of probability measure sequence is equivalent to the convergence of its Stieltjes transform sequence towards the corresponding transform of the limiting measure (for example, see \cite{BaiS}). So under the assumption $\mathbb{E}\{X_{11}^4\}<\infty$, one has when $N\rightarrow\infty$, $m_N(z)$ almost surely converges to the Stieltjes transform $f(z)$ of the semicircle law given by
\begin{eqnarray}
f(z)=\frac{-z+\sqrt{z^2-4\sigma^4}}{2\sigma^2}. \label{1.5}
\end{eqnarray}
Here the square root is specified to be the one with positive imaginary part. What's more, due to the basic fact $|m_N(z)|\leq\eta^{-1}$, the a.s convergence of $m_N(z)$ also implies
\begin{eqnarray}
f_N(z):=\mathbb{E}\{m_N(z)\}=f(z)+o(1).\label{1.0}
\end{eqnarray}

In this paper, we will take a step further to calculate the leading order term of the remainder. A corresponding work for the Wigner matrices has been taken in \cite{KKP} by Khorunzhy et al., and extended by some subsequent articles to some other matrix ensembles, see\cite{APS}, \cite{KK}. Such a topic is often referred to as  the asymptotic expansion of $f_N(z)$ in the literature. For ease of the expression and calculation, we just consider the normalized case of $\sigma^2=1$ hereafter,  so the semicircle law in (\ref{1.2}) will be a standard one, the general case is just analogous. We denote the $\alpha$-th moment of $X_{11}$ by $\omega_\alpha$ and $\alpha$-th cumulant by $\kappa_\alpha$, particularly there are $\kappa_2=\sigma^2=1$ and $\kappa_4=\omega_4-3$ in our setting. Our first main theorem is on the asymptotic expansion of $f_N(z)$:\\\\
 $\bf{Theorem~1.1}$
 \emph{Consider the matrix model $H_N$ defined in (\ref{1.1}) with $\sigma^2=1$. And $X=[X_{jk}]_{M\times N}$ is an $M\times N$ real random matrix with mean zero and variance $1$ i.i.d entries. We assume $M=M(N)$ and $N/M\rightarrow 0$ as $N\rightarrow \infty$. If  $\mathbb{E}|X_{11}|^{5+\delta}<\infty$ with some small positive $\delta$ , then we have the following asymptotic expansion holds for any fixed $z$ with $\eta>0$
 \begin{eqnarray*}
 f_N(z)&=&f(z)\bigg{\{}1-\sqrt{\frac NM}\frac{f^3(z)}{1-f^2(z)}+\frac1N\bigg{(}\frac{f^2(z)}{(1-f^2(z))^2}\nonumber\\
 &+&\kappa_4\frac{f^2(z)}{1-f^2(z)}\bigg{)}\bigg{\}}+o\left(\sqrt{\frac NM}\vee\frac1N\right). \label{1.6}
 \end{eqnarray*}}\\
$\emph{Remark~1.1}$:  By comparing with the corresponding result for the Wigner matrix in \cite{KKP}, we can find a different first order remainder with order $\sqrt{\frac NM}$ rather than $1/N$ appears when $M=o(N^3)$. What's more, we do only require $\eta>0$ in our result, rather than $\eta\geq 2\omega$ for Wigner matrix case in \cite{KKP}.

 (\ref{1.5}) and (\ref{1.6}) show the convergence of $m_N(z)$ a.s. and in expectation respectively, as usual, we shall consider the fluctuation of $m_N(z)$ as a further job. In this paper, we will deal with the more general linear eigenvalue statistic
 \begin{eqnarray*}
 \mathcal{L}_N[\varphi]=\sum_{i=1}^{N}\varphi(\lambda_i) \label{1.7}
 \end{eqnarray*}
 of $H_N$ with the test function $\varphi$ satisfies some smooth conditions. Note that $Nm_N(z)$ is a linear eigenvalue statistic with $\varphi(x)=1/(x-z)$. The topic of CLT for linear eigenvalue statistics is really classical and attractive in the Random Matrix Theory, there are a vast of related articles for different matrix ensembles, for examples, see \cite{AZ}, \cite{BaS}, \cite{Joh}, \cite{LP}, \cite{SM}, \cite{SS}. However, for the large dimensional sample covariance matrices, there is no corresponding result on the case of $N/M\rightarrow0$. So as a complement, we will discuss it in this paper. For convenience, we follow the notations in \cite{SM} to denote $\xi^\circ=\xi-\mathbb{E}\{\xi\}$ for any random variable $\xi$ below. And introduce the norm
\begin{eqnarray*}
||\varphi||_s^2=\int(1+2|k|)^{2s}|\widehat{\varphi}(k)|^2dk,\qquad \widehat{\varphi}(k)=\frac{1}{2\pi}\int e^{ikx}\varphi(x)dx. \label{1.8}
\end{eqnarray*}
Then we can state our second result as follows:\\\\
 $\bf{Theorem~1.2}$\emph{ Consider the matrix model $H_N$ defined in (\ref{1.1}) with $\sigma^2=1$. And $X=[X_{jk}]_{M\times N}$ is an $M\times N$ real random matrix with mean zero and variance $1$ i.i.d entries. We assume $M=M(N)$ and $N/M\rightarrow 0$ as $N\rightarrow \infty$. If $\mathbb{E}|X_{11}|^{5}<\infty$, then for any real test function $\varphi$ satisfying $||\varphi||_{3/2+\epsilon}<\infty$ with any $\epsilon>0$, we have $\mathcal{L}_N^\circ[\varphi]$ converges weakly to the centered Gaussian distribution with variance
 \begin{eqnarray}
 V[\varphi]&=&\frac{1}{2\pi^2}\int_{-2}^{2}\int_{-2}^{2}\left(\frac{\varphi(\lambda_1)-\varphi(\lambda_2)}{\lambda_1-\lambda_2}\right)^2
 \frac{4-\lambda_1\lambda_2}{\sqrt{4-\lambda_1^2}\sqrt{4-\lambda_2^2}}d\lambda_1d\lambda_2\nonumber\\
 &&+\frac{\kappa_4}{4\pi^2}\left(\int_{-2}^{2}\frac{\varphi(\mu)\mu}{\sqrt{4-\mu^2}}d\mu\right)^2. \label{1.9}
 \end{eqnarray}}
 $\emph{Remark~1.2}$: Again we can compare it to the corresponding result of the Wigner matrices in \cite{SM}, the variance coincides with the counterpart for a Gaussian Wigner matrix (Not GOE!) with variance $\omega_4-1$ diagonal entries. This can be explained as a result of the central limit theorem and asymptotic independence effect for the matrix entries. We require $\varphi$ to be real in the above Theorem, so when we deal with $m_N(z)$, we need to work with $\mathrm{Re}m_N(z)$ and $\mathrm{Im}m_N(z)$ as in \cite{KKP}.

 Our article is organized as follows. We will provide some basic facts and tools in Section 2, and as a warm-up, we use them to revisit the semicircle law for Gaussian case (i.e. $X_{11}$ is Gaussian) in an average sense. In Section 3 we present the proof of Theorem$~1.1$, which is based on the preliminaries provided in Section 2 together with some extra lemmas, whose proofs are postponed to section 4. In Section 5, we present the proof of Theorem$~1.2$ in short and leave all technical details in the Appendix. Throughout the paper, we will use $C$ and $C_i(i=1,2,3,4)$ to denote positive constants which may be different from line to line, sometimes may depend on $\eta$.\\

 \section{Preliminaries And Gaussian Case}
In this section we present some basic tools and facts required in the sequel, and use them to revisit the semicircle law proved in \cite{BY} in an average sense. For convenience, when there is no confusion, we will get rid of the subscript $N$ in the notations of matrices and $z$ as a variable in the notation $G_N(z)$ below. If we set $Y=[Y_{jk}]_{M\times N}:=(MN)^{-1/4}X$, we have
\begin{eqnarray*}
H&=&Y^TY-\sqrt{\frac MN}I_N. \label{2.1}
\end{eqnarray*}
To do the job we state the following two lemmas without detailed but routine proofs.\\\\
$\bf{Lemma~2.1}$ \emph{(Generalized Stein's equation)\\
For any real-valued random variable $\xi$ with $\mathbb{E}\{|\xi|^{p+2}\}<\infty$ and complex valued function $g(t)$ with continuous and bounded $p+1$ derivatives, we have
\begin{eqnarray}
 \mathbb{E}\{\xi g(\xi)\}=\sum_{a=0}^{p}\frac{\kappa_{a+1}}{a!}\mathbb{E}\{g^{(a)}(\xi)\}+\epsilon, \label{2.2}
 \end{eqnarray}
where $\kappa_a$ is the $a$-th cumulant of $\xi$, and
 \begin{eqnarray*}
 |\epsilon|\leq C\sup_{t}|g^{(p+1)}(t)|\mathbb{E}\{|\xi|^{p+2}\}, \label{2.3}
 \end{eqnarray*}
where the positive constant $C$ depends on $p$.}
\begin{rem}
The proof of Lemma~$2.1$ can be found in a lot of references, for example, see \cite{KKP} for details.
\end{rem}
When $\xi$ is centered Gaussian, (\ref{2.2}) reduces to the famous Stein's equation:

\begin{eqnarray}
\mathbb{E}\{\xi g(\xi)\}=\mathbb{E}\xi^2\cdot\mathbb{E}\{g'(\xi)\}. \label{2.4}
\end{eqnarray}

Our second lemma is on the derivatives $D_{jk}\{\cdot\}$ with respect to the matrix entry $Y_{jk}$, which will be used frequently in Section 2 and 3.\\\\
$\bf{Lemma~2.2}$ \emph{For any $\alpha, j\in \{1,2,\cdots,M \}$ and $\beta,k\in\{1,2,\cdots,N\}$, we have\\
(i):
\begin{eqnarray*}
D_{jk}G_{\alpha\beta}=-(YG)_{j\alpha}G_{\beta k}-(YG)_{j\beta}G_{\alpha k}, \label{2.5}
\end{eqnarray*}
(ii):
\begin{eqnarray*}
D_{jk}(YG)_{\alpha\beta}=\delta_{\alpha j}G_{\beta k}-G_{\beta k}(YGY^T)_{j\alpha}-(YG)_{j\beta}(YG)_{\alpha k},\label{2.6}
\end{eqnarray*}
(iii):
\begin{eqnarray*}
D_{jk}[G_{kk}(YGY^T)_{jj}]=2G_{kk}(YG)_{jk}-4G_{kk}(YG)_{jk}(YGY^T)_{jj},\label{2.7}
\end{eqnarray*}
where $\delta_{\alpha j}$ in (ii) is the Kronecker delta function.}
\begin{rem}
The proof of Lemma~$2.2$ is based on the following resolvent identity for real symmetric matrix $A$ and $B$
\begin{eqnarray*}
(A+B-z)^{-1}-(A-z)^{-1}=-(A+B-z)^{-1}B(A-z)^{-1},
\end{eqnarray*}
which implies
\begin{eqnarray}
\frac{d}{d\epsilon}(A+\epsilon B-z)^{-1}|_{\epsilon=0}=-(A-z)^{-1}B(A-z)^{-1}. \label{2.15}
\end{eqnarray}
 Now we set $A=H$ and $A+\epsilon B=(Y+\epsilon E(j,k))^T(Y+\epsilon E(j,k))$. Here $E(j,k)$ represents the $M\times N$ matrix with $(j,k)$-th entry to be $1$ and others $0$. Then by (\ref{2.15}) it is not difficult to get $(i)$ of Lemma$~2.2$. And $(ii)$ and $(iii)$ can be proved by the chain rule. We omit the details here, in fact it is quite similar with the counterparts in \cite{KKP} and \cite{LP}.
\end{rem}

As a warm-up for the main task in Section 3, we use the above two lemmas to prove $f_N(z)$ tends to $f(z)$ (as $N\rightarrow\infty$) for the Gaussian case below. By the basic relation between a matrix and its Green function
\begin{eqnarray*}
G=-z^{-1}+z^{-1}GH,\label{2.8}
\end{eqnarray*}
we have
\begin{eqnarray}
f_N(z)&=&-z^{-1}+z^{-1}\frac1N\mathbb{E}\{TrGH\}\nonumber\\
&=&-z^{-1}-z^{-1}\sqrt{\frac MN}f_N(z)+z^{-1}\frac1N\sum_{j,k}\mathbb{E}\{Y_{jk}(YG)_{jk}\}.\label{2.9}
\end{eqnarray}
Then we can use the Stein's equation (\ref{2.4}) to (\ref{2.9}), which yields
\begin{eqnarray}
f_N(z)&=&-z^{-1}-z^{-1}\sqrt{\frac MN}f_N(z)\nonumber\\
      &+&z^{-1}\frac1N\frac{1}{\sqrt{MN}}\sum_{j,k}\mathbb{E}\{G_{kk}-((YG)_{jk})^2-G_{kk}(YGY^T)_{jj}\}\nonumber\\
      &=&-z^{-1}-z^{-1}\frac1N\frac{1}{\sqrt{MN}}\mathbb{E}\{TrYG^2Y^T\}-z^{-1}\frac{1}{\sqrt{MN}}\mathbb{E}\{m_N(z)TrYGY^T\}\nonumber\\
      &=&-z^{-1}-z^{-1}\frac{N+1}{\sqrt{MN}}f_N(z)-z^{-1}(1+\sqrt{\frac NM}z)\mathbb{E}\{m_{N}^{2}(z)\}\nonumber\\
      &~&-z^{-1}\frac{1}{N^2}(1+\sqrt{\frac NM}z)\mathbb{E}\{TrG^2\},\label{2.10}
\end{eqnarray}
where we have used the fact that
\begin{eqnarray*}
&&TrYGY^T=Tr\bigg{(}I_N+(\sqrt{\frac MN}+z)G\bigg{)},\nonumber\\
&&TrYG^2Y^T=Tr\bigg{(}G+(\sqrt{\frac MN}+z)G^2\bigg{)}.
\end{eqnarray*}
Note the trivial bound $||G||\leq \eta^{-1}$ for the matrix norm of $G$, (\ref{2.10}) implies
\begin{eqnarray*}
f_N(z)=-z^{-1}-z^{-1}\mathbb{E}\{m_{N}^{2}(z)\}+O\bigg{(}\frac1N\vee\sqrt{\frac NM}\bigg{)}.\label{2.11}
\end{eqnarray*}

To estimate $\mathbb{E}\{m_{N}^{2}(z)\}$, we need to derive a bound for $\mathrm{Var}\{m_N(z)\}$.
By using the Poincar\'{e} inequality for Gaussian matrix entries as in Proposition$~2.4$ of \cite{LP}, it is not difficult to get  $$\mathrm{Var}\{m_N(z)\}=O\bigg{(}\frac{1} {N^2}\bigg{)}.$$
Consequently one has
\begin{eqnarray}
f_N(z)=-z^{-1}-z^{-1}f_N^2(z)+O\bigg{(}\frac1N\vee\sqrt{\frac NM}\bigg{)}.\label{2.12}
\end{eqnarray}
On the other hand, it is well known that $f(z)$ satisfies the equation
\begin{eqnarray}
f(z)=-z^{-1}-z^{-1}f^2(z).\label{2.13}
\end{eqnarray}
By a routine comparison issue on (\ref{2.12}) and (\ref{2.13}) we have $f_N(z)\rightarrow f(z)$ as $N\rightarrow\infty$ (see Section 2 of \cite{BaiS} for example), which implies the convergence of the expected empirical spectral distribution to the semicircle law.

However, for Theorem$~1.1$, we need to take a step further to calculate the leading order term of the remainder $f_N(z)-f(z)$ precisely for the general distribution case. So at first we need use (\ref{2.2}) instead of (\ref{2.4}). As a result, more involved estimate towards the derivatives is required. What's more, the Poincar\'{e} inequality is no longer valid, we need to estimate $\mathrm{Var}\{m_N(z)\}$ by a martingale difference method as in \cite{SM}.
\section{Asymptotic Expansion For $f_N(z)$}
To prove Theorem$~1.1$, we begin with the basic idea raised in \cite{KKP}, use (\ref{2.2}) to (\ref{2.9}) and estimate the derivatives. However, the standard process will be too involved for our matrix model $H_N$. As is shown in \cite{KKP}, for the Wigner matrix, every term in the expansion formula can be factorized into the entries of $G_N(z)$ which can be bounded by $\eta^{-1}$. Owing to such a trivial bound, many estimates hold obviously, especially for the remainder term. And the approach of iterating the expansion process also works well. However as we will see, for our case, the derivatives are in more complicated forms. There are no trivial bound for some factors of the terms (see (\ref{3.13}), (\ref{3.14})), as well, iterating the expansion process will bring factors with new types. And to bound the remainder, we need to use a truncation technic at first.

Now we truncate $X_{jk}$ at $\tau:=(MN)^{1/4-t}$ with some small positive $t\leq\frac{\delta}{100}$ (say). Then we denote the empirical spectral distribution of the matrix defined in (\ref{1.1}) with truncated entries $X^\tau_{jk}:=X_{jk}\mathbf{1}_{\{|X_{jk}|\leq \tau\}}$ by $F^\tau_N(x)$. Further we denote the stieltjes transform of $F^\tau_N(x)$ by $m^\tau_N(z)$ and set $f^\tau_N(z)=\mathbb{E}m^\tau_N(z)$.
By taking into account that $|m_N(z)|,|m^\tau_N(z)|\leq \eta^{-1}$ and
\begin{eqnarray*}
\sum_{jk}\mathbb{P}\{|X_{jk}|\geq \tau\}\leq CMN\cdot \tau^{-(5+\delta)}\leq C(MN)^{-\frac 14-t}=o\bigg{(}\sqrt{\frac NM}\vee\frac1N\bigg{)},
\end{eqnarray*}
we have
$$
|f_N(z)-f^\tau_N(z)|=o\bigg{(}\sqrt{\frac NM}\vee\frac1N\bigg{)}.
$$
 After the truncation, all of the first 5 moments of $X_{jk}$ will be modified. For example,
  \begin{eqnarray*}
  |\mathbb{E}\{X_{jk}\mathbf{1}_{\{|X_{jk}|\leq \tau\}}\}|=|\mathbb{E}\{X_{jk}\mathbf{1}_{\{|X_{jk}|\geq \tau\}}\}| \leq\tau^{-(4+\delta)}\mathbb{E}|X_{jk}|^{5+\delta}
 \leq C(MN)^{-1-t}.
  \end{eqnarray*}

  Similarly we can calculate the modification of the $\alpha$-th moments of $X_{jk}$ under the truncation for all $\alpha=2,3,4,5$. So without loss of generality, we can always assume $X_{j,k}, (j,k)\in\{1,\cdots,M\}\times\{1,\cdots,N\}$ are i.i.d and
 \begin{eqnarray}
 &&|X_{jk}|\leq (MN)^{1/4-t}, \quad |\mathbb{E}\{X_{jk}\}|\leq C(MN)^{-1-t},\nonumber\\
 &&\mathbb{E}\{X_{jk}^2\}=1+C(MN)^{-\frac 34-t},
 \quad \mathbb{E}\{X_{jk}^\alpha\}=\omega_\alpha+o(1), \alpha=3,4,5 \nonumber\\
 &&\mathbb{E}\{|X_{jk}|^{5+\delta}\}\leq C.\label{3.1.4}
 \end{eqnarray}
 And for simplicity of the notations, we still use $f_N(z)$ to denote $f^\tau_N(z)$ below.
 What's more, we need to present two lemmas before the rigorous proof of Theorem$~1.1$. The first lemma is on the estimate of $\mathrm{Var}\{m_N(z)\}$, as we have seen in Section 2, it is a necessary ingredient for the final result. The second lemma is on an estimate towards every diagonal entry of the Green function, which can help to overcome effectively the difficulties we have mentioned above. Their proofs will be postponed to Section 4.\\\\
$\bf{Lemma~3.1}$: \emph{Under the assumptions of $\omega_1=0$ and $\omega_4<\infty$, we have
\begin{eqnarray*}
\mathrm{Var}\{m_N(z)\}=O\bigg{(}\frac{1}{N^2}\bigg{)} \label{3.1}
\end{eqnarray*}
for any fixed $z$ with $\eta>0$}.\\\\
To present the second lemma, we need to introduce the matrix
\begin{eqnarray*}
\widetilde{G}=[\widetilde{G}_{jk}]_{M\times M}:=\bigg{(}YY^T-\sqrt{\frac MN}-z\bigg{)}^{-1}. \label{3.2}
\end{eqnarray*}\\
$\bf{Lemma~3.2}$: \emph{Under the assumptions of Lemma$~3.1$, for any fixed $z$ with $\eta>0$ we have the following two estimates\\
(i): for any $k\in\{1,2,\cdots, N\}$
\begin{eqnarray}
\mathbb{E}\left|G_{kk}+\frac{1}{z+f_N(z)}\right|^2\leq C_1\frac 1N+C_2\frac NM \label{3.3}
\end{eqnarray}
(ii): for any $j\in\{1,2,\cdots, M\}$
\begin{eqnarray}
\mathbb{E}\bigg{|}\widetilde{G}_{jj}+\frac{1}{\sqrt{\frac MN}+z+f_N(z)}\bigg{|}^2\leq C_1(\frac NM)^3+C_2\frac{N}{M^2} \label{3.4}
\end{eqnarray}}
$\emph{Remark~3.1}$ In fact, we can replace $(z+f_N(z))^{-1}$ in (\ref{3.3}) by $f_N(z)$ or even $f(z)$ with the same bound on the right hand side, but as a lemma, we shall not do that before we get the final conclusion on the differences between $(z+f_N(z))^{-1}$, $f_N(z)$ and $f(z)$.  \\
\begin{proof}[\emph{Proof of Theorem$~1.1$}]
Use (\ref{2.2}) to (\ref{2.9}), we have the following expansion for the general distribution case
\begin{eqnarray}
f_N(z)&=&-z^{-1}-z^{-1}\sqrt{\frac{M}{N}}f_N(z)\nonumber\\
      &~&+z^{-1}\frac{1}{N}\sum_{a=0}^{3}\frac{1}{(MN)^{\frac{a+1}{4}}}\sum_{j,k}\frac{\kappa_{a+1}}{a!}
      \mathbb{E}{D_{jk}^{a}(YG)_{jk}}+R_N, \label{3.5}
\end{eqnarray}
where $D_{\alpha\beta}^{a}\{\cdot\}$ represents the $a$-th derivative with respect to $Y_{\alpha\beta}$ and
\begin{eqnarray*}
R_N\leq C\frac1N\frac{1}{(MN)^{\frac54}}\sum_{j,k}\sup_{jk}\mathbb{E}_{jk}|D^4_{jk}(YG)_{jk}|. \label{3.6}
\end{eqnarray*}
Observe that in (\ref{3.5}), the $a=0$ term
\begin{eqnarray*}
z^{-1}\frac{1}{N}\frac{1}{(MN)^{1/4}}\mathbb{E}\{X_{11}\}\sum_{jk}\mathbb{E}(YG)_{jk}=o\bigg{(}\sqrt{\frac{M}{N}}\vee\frac{1}{N}\bigg{)}
\end{eqnarray*}
by taking into account (\ref{3.1.4}) and the estimation (\ref{3.17.5}) proved below. So we will focus on the estimations of $a\geq 1$ terms and $R_N$ in the sequel.
Here $\sup_{jk}$ means the supremum is taken w.r.t $Y_{jk}$ and $E_{jk}$ represents the conditional expectation w.r.t $Y_{jk}$. By using Lemma$~2.2$ repeatedly we have
\begin{eqnarray}
D_{jk}(YG)_{jk}&=&G_{kk}-[(YG)_{jk}]^2-G_{kk}(YGY^T)_{jj}\nonumber\\
D_{jk}^2(YG)_{jk}&=&-6G_{kk}(YG)_{jk}+6G_{kk}(YG)_{jk}(YGY^T)_{jj}+2[(YG)_{jk}]^3\label{3.8}\\
D_{jk}^3(YG)_{jk}&=&-6(G_{kk})^2+36G_{kk}[(YG)_{jk}]^2\nonumber\\
&~&+12(G_{kk})^2(YGY^T)_{jj}-36G_{kk}[(YG)_{jk}]^2(YGY^T)_{jj}\nonumber\\
&~&-6(G_{kk})^2[(YGY^T)_{jj}]^2-6[(YG)_{jk}]^4\label{3.9}\\
D_{jk}^4(YG)_{jk}&=&120[G_{kk}]^2(YG)_{jk}-240G_{kk}[(YG)_{jk}]^3\nonumber\\
&~&-240[G_{kk}]^2(YG)_{jk}(YGY^T)_{jj}+240G_{kk}[(YG)_{jk}]^3(YGY^T)_{jj}\nonumber\\
&~&+120[G_{kk}]^2(YG)_{jk}[(YGY^T)_{jj}]^2+24[(YG)_{jk}]^5\nonumber
\end{eqnarray}
for any $(j,k)\in\{1,\cdots,M\}\times\{1,\cdots,N\}$. Observe that for any integer $a\geq 1$, $D_{jk}^a(YG)_{jk}$ is a linear combination of terms in the form of
\begin{eqnarray}
F_{jk}(a_1,a_2,a_3):=[G_{kk}]^{a_1}[(YG)_{jk}]^{a_2}[(YGY^T)_{jj}]^{a_3}\label{3.11}
\end{eqnarray}
with nonnegative integers $a_1, a_2, a_3$ satisfying $a_1+a_2+a_3\leq a+1$ . Moreover, by Lemma$~2.2$,
it is easy to see one $(YGY^T)_{jj}$ must appears together with one $G_{kk}$, thus we always have $a_1\geq a_3$.
By observing that
\begin{eqnarray*}
YGY^T=\widetilde{G}YY^T=I_M+(\sqrt{\frac{M}{N}}+z)\widetilde{G}\label{3.12}
\end{eqnarray*}
and the trivial fact $||\widetilde{G}||\leq \eta^{-1}$ we have
\begin{eqnarray}
|(YGY^T)_{jj}|\leq C\sqrt{\frac MN}.\label{3.13}
\end{eqnarray}
Similarly we have the following estimate
\begin{eqnarray}
|(YG)_{jk}|\leq\left(\sum_{k}|(YG)_{jk}|^2\right)^{1/2}=|(YGG^*Y^T)_{jj}|^{1/2}\leq C\left(\frac MN\right)^{1/4}.\label{3.14}
\end{eqnarray}
But these two bounds are too bad, thus we need to figure out the higher order term $R_N$ more carefully. Fortunately, when we take expectations in (\ref{3.5}), we can use the following estimate
\begin{eqnarray}
\mathbb{E}\bigg{|}(YGY^T)_{jj}-\frac{f_N(z)}{\sqrt{\frac MN}+z+f_N(z)}\bigg{|}^2\leq C_1(\frac NM)^2+C_2\frac{1}{M},\label{3.15}
\end{eqnarray}
which is a direct consequence of (\ref{3.4}) and (\ref{3.12}).

To deal with $R_N$, we observe that
 \begin{eqnarray}
| F_{jk}(a_1,a_2,a_3)|
&\leq& C|(YG)_{jk}|^{a_2}|(YGY^T)_{jj}|^{a_3}\nonumber\\
&=&C|\sum_{\alpha}Y_{j\alpha}G_{\alpha k}|^{a_2}|\sum_{\beta,\gamma}Y_{j\beta}Y_{j\gamma}G_{\beta\gamma}|^{a_3}\nonumber\\
&\leq&\bigg{(}(\sum_{\alpha} Y_{j\alpha}^2)(\sum_{\alpha}|G_{\alpha j}|^2)\bigg{)}^{\frac{a_2}{2}}\bigg{(}(\sum_{\beta,\gamma}Y_{j\beta}^2Y_{j\gamma}^2)
(\sum_{\beta,\gamma}|G_{\beta\gamma}|^2)\bigg{)}^{\frac{a_3}{2}}\nonumber\\
&\leq&CN^{\frac{a_3}{2}}(\sum_{\alpha}Y_{j\alpha}^2)^{\frac{a_2+2a_3}{2}}.\label{3.12.4}
 \end{eqnarray}
 In the second inequality we used Cauchy Schwartz inequality and in the last step we used the fact that
 \begin{eqnarray*}
 &&\sum_{\alpha}|G_{\alpha j}|^2=[G^*G]_{jj}\leq C,\quad \sum_{\beta,\gamma}|G_{\beta\gamma}|^2=TrG^*G=CN\\
 &&\sum_{\beta,\gamma}Y_{j\beta}^2Y_{j\gamma}^2=\left(\sum_{\alpha}Y_{j\alpha}^2\right)^2.
 \end{eqnarray*}
 So by (\ref{3.12.4}) we have
\begin{eqnarray*}
\sup_{jk}\mathbb{E}_{jk}| F_{jk}(a_1,a_2,a_3)|
&\leq& CN^{\frac{a_3}{2}}\bigg{(}\sup|Y_{jk}|^{a_2+2a_3}+\mathbb{E}\bigg{(}\sum_{\alpha\neq k}Y_{j\alpha}^2\bigg{)}^{\frac{a_2+2a_3}{2}}\bigg{)}\nonumber\\
&\leq& CN^{\frac{a_3}{2}}\cdot(MN)^{-t(a_2+2a_3)}+CN^{\frac{a_3}{2}}\cdot N^{\frac{a_2+2a_3}{2}}(MN)^{-\frac{a_2+2a_3}{4}}
\end{eqnarray*}
Then it is not difficult to see
$$|R_N|=o(\sqrt{\frac NM}\vee\frac1N)$$
by taking into account $a_2+2a_3\leq 5$.

Now we can turn our attention to (\ref{3.8}) and (\ref{3.9}). Firstly, we shall provide some coarse estimate on $\sum_{jk}|(YG)_{jk}|^{a_2}$. When $a_{2}\geq 2$ we have
\begin{eqnarray}
|\sum_{j,k}[(YG)_{jk}]^{a_2}|&\leq& \sum_{j,k}(|(YG)_{jk}|^2)^{\frac{a_2}{2}}\nonumber\\
                       &=&\sum_{j,k}[(G^*Y^T)_{kj}(YG)_{jk}]^{\frac{a_2}{2}}\nonumber\\
                       &\leq&\sum_{k}[(G^*Y^TYG)_{kk}]^{\frac{a_2}{2}}\nonumber\\
                       &=&\sum_{k}\bigg{[}G_{kk}+(\sqrt{\frac{M}{N}}+\bar{z})(G^*G)_{kk}\bigg{]}^{\frac{a_2}{2}}\nonumber\\
                       &=&O(M^{\frac{a_2}{4}}N^{1-\frac{a_2}{4}}).\label{3.16}
\end{eqnarray}
When $a_{2}=1$, we can use the elementary inequality
\begin{eqnarray}
|\sum_{j,k}(YG)_{jk}|\leq\bigg{(}MN\sum_{jk}|(YG)_{jk}|^2\bigg{)}^{\frac12}=O\bigg{(}(MN)^{\frac34}\bigg{)}.\label{3.17.5}
\end{eqnarray}
And for this time we only need to bound the expectations, so (\ref{3.15}) can be used. By Lemma$~3.2$, (\ref{3.15}),(\ref{3.16}) and (\ref{3.17.5}), we can discard all the terms except for
\begin{eqnarray*}
-6F_{jk}(1,1,0),\qquad -6F_{jk}(1,0,0)\label{3.19}
\end{eqnarray*}
 from (\ref{3.8}), (\ref{3.9}) respectively. For $F_{jk}(1,1,0)$, according to Lemma$~3.2$ and (\ref{3.17.5}), it suffices to estimate $\sum_{jk}\mathbb{E}\{(YG)_{jk}\}$ instead. To get a lower bound than the coarse one (\ref{3.17.5}), we need to iterate the expansion process, the main term of $\sum_{jk}\mathbb{E}\{(YG)_{jk}\}$ is enough for us. We do it as follows:
\begin{eqnarray}
\mathbb{E}\sum_{j,k}(YG)_{jk}
&=&\sum_{j,\alpha,k}\mathbb{E}(Y_{j\alpha}G_{\alpha k})\nonumber\\
&=&\sum_{j,\alpha,k}\bigg{[}\frac{1}{\sqrt{MN}}\mathbb{E}D_{j\alpha}G_{\alpha k}+\frac{\kappa_3}{2}\frac{1}{(MN)^{\frac{3}{4}}}\mathbb{E}D_{j\alpha}^{2}G_{\alpha k}\bigg{]}+\epsilon_{N},\label{3.20}
\end{eqnarray}
where
\begin{eqnarray*}
|\epsilon_{N}|\leq C\frac{\omega_4}{MN}\sum_{j,\alpha,k}\sup_{j\alpha}\mathbb{E}_{j\alpha}|D_{j\alpha}^3G_{\alpha k}|.\label{3.21}
\end{eqnarray*}
By using Lemma$~2.2$ again we can get
\begin{eqnarray}
D_{j\alpha}G_{\alpha k}&=&-(YG)_{j\alpha}G_{\alpha k}-(YG)_{jk}G_{\alpha\alpha}\label{3.22}\\
D_{j\alpha}^2G_{\alpha k}&=&-2G_{\alpha\alpha}G_{\alpha k}+2G_{\alpha k}[(YG)_{j\alpha}]^2\nonumber\\
&~&+2G_{\alpha\alpha}G_{\alpha k}(YGY^T)_{jj}+4G_{\alpha\alpha}(YG)_{j\alpha}(YG)_{jk}\label{3.23}\\
D_{j\alpha}^3G_{\alpha k}&=&18G_{\alpha\alpha}G_{\alpha k}(YG)_{j\alpha}+6G_{\alpha\alpha}^2(YG)_{jk}\nonumber\\
&&-18G_{\alpha\alpha}G_{\alpha k}(YG)_{j\alpha}(YGY^T)_{jj}-6G_{\alpha k}[(YG)_{j\alpha}]^3\nonumber\\
&&-18G_{\alpha\alpha}(YG)_{jk}[(YG)_{j\alpha}]^2-6[G_{\alpha\alpha}]^2(YG)_{jk}(YGY^T)_{jj}.\nonumber
\end{eqnarray}
Inserting (\ref{3.22}) and (\ref{3.23}) into (\ref{3.20}) one finds that
\begin{eqnarray*}
\mathbb{E}\sum_{j,k}(YG)_{jk}&=&-\frac{1}{\sqrt{MN}}\sum_{j,k}\mathbb{E}\bigg{[}(YG^2)_{jk}+(YG)_{jk}TrG\bigg{]}\\
&~&~~+\frac{\kappa_{3}}{(MN)^{\frac34}}\sum_{\alpha,k}\mathbb{E}\bigg{[}-MG_{\alpha\alpha}G_{\alpha k}+G_{\alpha k}(GY^TYG)_{\alpha\alpha}\\
&~&~~G_{\alpha\alpha}G_{\alpha k}Tr(YGY^T)+2(GY^TYG)_{\alpha k}G_{\alpha\alpha}\bigg{]}
+\epsilon_{N}.
\end{eqnarray*}
Similarly to (\ref{3.17.5}) we can provide
\begin{eqnarray*}
\mathbb{E}\sum_{j,k}(YG^2)_{jk}=O\left((MN)^{\frac{3}{4}}\right).\label{3.25}
\end{eqnarray*}
The other terms in the expansion are easier to estimate by a similar calculation as we have done above. Moreover, similar to the case of $R_N$ it is not difficult to check that $\epsilon_N$ has no contribution to the main term of $\sum_{j,k}\mathbb{E}\{(YG)_{jk}\}$. So we can get a bound as
\begin{eqnarray*}
\mathbb{E}\sum_{j,k}(YG)_{jk}=O\left(M^{\frac{1}{4}}N^{\frac{5}{4}}\right).\label{3.26}
\end{eqnarray*}

Therefore, under the estimates above, we arrive at
\begin{eqnarray*}
f_N(z)&=&-z^{-1}-z^{-1}\frac1N\mathbb{E}\sum_{j,k}[(YG)_{jk}]^2\nonumber\\
&&-z^{-1}\frac1N\mathbb{E}\sum_{j,k}G_{kk}(YGY^T)_{jj}-z^{-1}\kappa_4\frac{1}{N^2}\mathbb{E}\sum_{k}G_{kk}^2+o(\frac1N).\label{3.27}
\end{eqnarray*}
By (\ref{2.10}) and Lemma$~3.1$, we can rewrite it as
\begin{eqnarray*}
f_N(z)&=&-z^{-1}-z^{-1}\sqrt{\frac NM}f_N(z)-z^{-1}(1+\sqrt{\frac NM}z)f_N^2(z)\nonumber\\
      &~&-z^{-1}\frac{1}{N^2}\mathbb{E}\{TrG^2(z)\}\nonumber\\
&~&-z^{-1}\kappa_4\frac{1}{N^2}\mathbb{E}\sum_{k}G_{kk}^2+o(\frac{1}{N}),\label{3.28}
\end{eqnarray*}
which implies
\begin{eqnarray}
\bigg{|}f_N(z)+\frac{1}{z+f_N(z)}\bigg{|}=O \bigg{(}\frac1N\vee\sqrt{\frac NM}\bigg{)}. \label{3.31}
\end{eqnarray}

Furthermore, we know both $f_N(z)$ and $f(z)$ are analytic function of $z$. Then by (\ref{1.0}) and the Cauchy's formula, we also have
\begin{eqnarray*}
\frac{d}{dz}(f_N(z))=\frac{d}{dz}(f(z))+o(1).
\end{eqnarray*}
Now we note that
\begin{eqnarray*}
\frac1N\mathbb{E}TrG^2(z)=\frac{d}{dz}(f_N(z))=\frac{d}{dz}(f(z))+o(1)=-\frac{f(z)}{z+2f(z)}+o(1).\label{3.29}
\end{eqnarray*}
Together with (\ref{3.3}) we finally have
\begin{eqnarray*}
f_N(z)&=&f(z)-\sqrt{\frac{N}{M}}\frac{(f(z)+zf^2(z))}{z+2f(z)}\nonumber\\
      &~&+\frac1N\left[\frac{f(z)}{(z+2f(z))^2}-\frac{\kappa_4f^2(z)}{z+2f(z)}\right]+O(\frac NM)+o\left(\frac1N\right).\label{3.30}
\end{eqnarray*}
This concludes the proof by inserting the basic relation $z+2f(z)=f(z)-1/f(z)$.
\end{proof}

\section{Variance Estimates}
In this section, we will state the proofs of Lemma$~3.1$ and $~3.2$. In fact, we can use the asymptotic expansion method again to estimate $\mathrm{Var}\{m_N(z)\}$, but it is really more complicated than in the Wigner case. We will use a martingale difference method used in \cite{SM} very recently. For convenience we introduce the notation $\gamma_N$ to replace $TrG$ below.
\begin{proof}[\emph{Proof of Lemma~3.1}]
If we denote $\mathbb{E}_{\leq k}\{\cdot\}$ and $\mathbb{E}_{k}\{\cdot\}$ as the expectation w.r.t the random variables of the first $k$ columns and $k$-th column of $Y_N$ respectively, by the classical martingale method in \cite{DFJ}, one has
\begin{eqnarray}
Var\{\gamma_N\}&=&\sum_{k=1}^{N}\mathbb{E}\{|\mathbb{E}_{\leq k-1}\{\gamma_N\}-\mathbb{E}_{\leq k}\{\gamma_N\}|^2\}\nonumber\\
             &=&\sum_{k=1}^{N}\mathbb{E}\{|\mathbb{E}_{\leq k-1}\{\gamma_N-\mathbb{E}_{k}\{\gamma_N\}\}|^2\}\nonumber\\
             &\leq&\sum_{k=1}^{N}\mathbb{E}\{|\gamma_N-\mathbb{E}_k\{\gamma_N\}|^2\}.\label{4.1}
\end{eqnarray}
Now we let $y_{k}$ be the $k$-th column of $Y_{N}$ and $B_{k}$ the $M\times(N-1)$ matrix consisting of the other $N-1$ columns of $Y_{N}$, then we have
 $$
H=\left(
\begin{array}{ccccc}
 y_{1}\cdot y_{1}-\sqrt{\frac{M}{N}} &(B_{1}^Ty_{1})^{T}  \\
 B_{1}^{T}y_{1} &B_{1}^{T}B_{1}-\sqrt{\frac{M}{N}}I_{N-1}  \\
\end{array}
\right).
$$
With the notations
$$H^{(k)}=B_{k}^{T}B_{k}-\sqrt{\frac{M}{N}}I_{N-1}$$
and $G^{(k)}=(H^{(k)}-z)^{-1}$, we have
\begin{eqnarray}
TrG-TrG^{(1)}=\frac{1+y_1\cdot B_1B_1^T\left(B_1B_1^T-\sqrt{\frac MN}-z\right)^{-2}y_1}{y_1\cdot y_1-\sqrt{\frac MN}-z-y_1\cdot B_1B_1^T\left(B_1B_1^T-\sqrt{\frac MN}-z\right)^{-1}y_1}=:\frac{1+U}{V}\label{4.2}
\end{eqnarray}
and $G_{11}=V^{-1}$. Here we used the basic identity
$$B_{1}\bigg{(}B_{1}^{T}B_{1}-\sqrt{\frac{M}{N}}-z\bigg{)}^{-n}B_{1}^{T}=B_{1}B_{1}^{T}\bigg{(}B_{1}B_{1}^{T}-\sqrt{\frac{M}{N}}-z\bigg{)}^{-n}, \qquad n=1,2.$$

To estimate (\ref{4.1}), we only need to deal with the first term
\begin{eqnarray}
\mathbb{E}\{|\gamma_N-\mathbb{E}_1\{\gamma_N\}|^2\}=\mathbb{E}\bigg{\{}\bigg{|}\frac{1+U}{V}-\mathbb{E}_1\bigg{\{}\frac{1+U}{V}\bigg{\}}\bigg{|}^2\bigg{\}},\label{4.3}
\end{eqnarray}
since the others are analogous. So it suffices to estimate $\mathbb{E}\{|UV^{-1}-\mathbb{E}_1\{UV^{-1}\}|^2\}$ and $\mathbb{E}\{|V^{-1}-\mathbb{E}_1\{V^{-1}\}|^2\}$. We will only present the estimate for the first one below. Clearly we have
\begin{eqnarray}
\mathbb{E}_1\{|UV^{-1}-\mathbb{E}_1\{UV^{-1}\}|^2\}&\leq& \mathbb{E}_1\{|UV^{-1}-\frac{\mathbb{E}_1\{U\}}{\mathbb{E}_1\{V\}}|^2\}\nonumber\\
&=&\mathbb{E}_1\bigg{\{}\bigg{|}\frac{U-\mathbb{E}_1\{U\}}{\mathbb{E}_1\{V\}}-\frac{V-\mathbb{E}_1\{V\}}{\mathbb{E}_1\{V\}}\frac{U}{V}\bigg{|}^2\bigg{\}}.\label{4.4}
\end{eqnarray}
By observing that $|\mathbb{E}_1\{V\}|^{-1}\leq \eta^{-1}$ and
\begin{eqnarray*}
|U/V|\leq|U/\mathrm{Im}V|=\frac{y_1\cdot B_1B_1^T|B_1B_1^T-\sqrt{\frac{M}{N}}-z|^{-2}y_1}{\eta+\eta y_1\cdot B_1B_1^T|B_1B_1^T-\sqrt{\frac{M}{N}}-z|^{-2}y_1}\leq\eta^{-1},
\end{eqnarray*}
 it suffices to provide the following estimates
\begin{eqnarray}
\mathbb{E}_1\{|U-\mathbb{E}_1\{U\}|^2\}, \quad
\mathbb{E}_1\{|V-\mathbb{E}_1\{V\}|^2\}
=O\left(\frac1N\right).\label{4.4.5}
\end{eqnarray}
For simplicity we introduce
\begin{eqnarray*}
M^{[1]}=[M^{[1]}_{jk}]_{M\times M}:=B_1B_1^T\left(B_1B_1^T-\sqrt{\frac MN}-z\right)^{-1}, \label{4.5}\\
M^{[2]}=[M^{[2]}_{jk}]_{M\times M}:=B_1B_1^T\left(B_1B_1^T-\sqrt{\frac MN}-z\right)^{-2}.\label{4.6}
\end{eqnarray*}
Thus we have
\begin{eqnarray*}
U-\mathbb{E}_1\{U\}&=&\sum_{i\neq j}M^{[2]}_{ij}Y_{i1}Y_{j1}+\sum_{i}M^{[2]}_{ii}(Y_{i1}^2-\mathbb{E}_1\{Y_{i1}^2\}),\label{4.7}\\
V-\mathbb{E}_1\{V\}&=&y_1\cdot y_1-\sqrt{\frac MN}-\sum_{i\neq j}M^{[1]}_{ij}Y_{i1}Y_{j1}-\sum_{i}M^{[1]}_{ii}(Y_{i1}^2-\mathbb{E}_1\{Y_{i1}^2\}).\label{4.8}
\end{eqnarray*}
It is not difficult to get
\begin{eqnarray}
&~&\mathbb{E}_1\{|U-\mathbb{E}_1\{U\}|^2\}
\leq C\frac{1}{MN}\sum_{i,j}|M^{[2]}_{ij}|^2=C\frac{1}{MN}Tr|M^{[2]}|^2,\label{4.9}\\
&~&\mathbb{E}_1\{|V-\mathbb{E}_1\{V\}|^2\}
\leq C_1\frac{1}{MN}Tr|M^{[1]}|^2+C_2\frac{1}{N}.\label{4.10}
\end{eqnarray}
Now we denote eigenvalues of $H^{(k)}$ by $\mu_{1}^{(k)} \cdots \mu_{N-1}^{(k)}$. Note that $\mu_{\alpha}^{(k)}(\alpha=1,\cdots,N-1)$ are also the eigenvalues of the matrix $$\check{H}^{(k)}:=B_{k}B_{k}^{T}-\sqrt{\frac{M}{N}}I_{M},$$
which has a $(M-N+1)$-multiple eigenvalue $-\sqrt{\frac{M}{N}}$. So it follows
\begin{eqnarray*}
Tr|M^{[n]}|^2=\sum_{\alpha=1}^{N-1}\frac{(\mu_{\alpha}^{(1)}+\sqrt{\frac MN})^2}{|\mu_{\alpha}^{(1)}-z|^{2n}}=O(M),\qquad n=1,2. \label{4.11}
\end{eqnarray*}
which implies (\ref{4.4.5}) by taking into account (\ref{4.9}) and (\ref{4.10}).
\end{proof}

\begin{proof}[\emph{Proof of Lemma$~3.2$}] We will only prove for $G_{11}$ and $\widetilde{G}_{11}$, the others are analogous. As we have shown above
\begin{eqnarray}
G_{11}=V^{-1}=\frac{1}{y_{1}\cdot y_{1}-\sqrt{\frac{M}{N}}-z-y_{1}\cdot M^{[1]}y_{1}}.\label{4.12}
\end{eqnarray}
We denote the unit eigenvector of $\check{H}^{(1)}$ corresponding to the eigenvalue $\mu_{\alpha}^{(1)}$ by $v_{\alpha}^{(1)}=(v_{\alpha}^{(1)}(1),\cdots,v_{\alpha}^{(1)}(M))$ $(\alpha=1,\cdots,N-1)$ and set
$$\xi_{\alpha}^{(1)}=\sqrt{MN}|y_{1}\cdot v_{\alpha}^{(1)}|^2.$$
It is easy to check $\mathbb{E}_1\{\xi_{\alpha}^{(1)}\}=1$. Then by (\ref{4.12}),
we have
\begin{eqnarray*}
G_{11}&=&\frac{1}{y_{1}\cdot y_{1}-\sqrt{\frac{M}{N}}-z-\frac{1}{\sqrt{MN}}\sum_{\alpha=1}^{N-1}\frac{(\mu_{\alpha}^{(1)}+\sqrt{\frac{M}{N}})
\xi_{\alpha}^{(1)}}{\mu_{\alpha}^{(1)}-z}}\nonumber\\
&=:&\frac{1}{-z-m_N(z)+r_1},\label{4.13}
\end{eqnarray*}
where
\begin{eqnarray*}
r_1&=&(y_1\cdot y_1-\sqrt{\frac{M}{N}})+\bigg{(}m_N(z)-\frac1N\sum_{\alpha=1}^{N-1}\frac{\sqrt{\frac NM}\mu_{\alpha}^{(1)}+1}{\mu_{\alpha}^{(1)}-z}\bigg{)}\nonumber\\
&~&-\frac1N\sum_{\alpha=1}^{N-1}\frac{\bigg{(}\sqrt{\frac NM}\mu_{\alpha}^{(1)}+1\bigg{)}(\xi_{\alpha}^{(1)}-1)}{\mu_{\alpha}^{(1)}-z}\nonumber\\
&=&\bigg{(}m_N(z)-\frac1N\sum_{\alpha=1}^{N-1}\frac{\sqrt{\frac NM}\mu_{\alpha}^{(1)}+1}{\mu_{\alpha}^{(1)}-z}\bigg{)}+V-\mathbb{E}_1\{V\}.\label{4.14}
\end{eqnarray*}
As a consequence we have
\begin{eqnarray*}
\left|G_{11}+\frac{1}{z+f_N(z)}\right|&=&\left|\frac{f_N(z)-m_N(z)+r_1}{(-z-m_N(z)+r_1)(z+f_N(z))}\right|\nonumber\\
\nonumber\\
&\leq& C|(f_N(z)-m_N(z))+r_1|.\label{4.15}
\end{eqnarray*}
Thus we obtain
\begin{eqnarray}
&~&\mathbb{E}\left|G_{11}+\frac{1}{z+f_N(z)}\right|^2\nonumber\\
&~&\leq C(\mathrm{Var}\{m_N(z)\}+\mathbb{E}\bigg{|}m_N(z)-\frac1N\sum_{\alpha=1}^{N-1}\frac{\sqrt{\frac NM}\mu_{\alpha}^{(1)}+1}{\mu_{\alpha}^{(1)}-z}\bigg{|}^2+\mathbb{E}|V-\mathbb{E}_1\{V\}|^2)\nonumber\\ \label{4.16}
\end{eqnarray}
Note the trivial bound
\begin{eqnarray}
\sqrt{\frac NM}\bigg{|}\frac1N\sum_{\alpha=1}^{N-1}\frac{\mu_{\alpha}^{(1)}}{\mu_{\alpha}^{(1)}-z}\bigg{|}=O(\sqrt{\frac NM}),\label{4.16.5}
\end{eqnarray}
and the fact that
\begin{eqnarray}
\bigg{|}m_{N}(z)-\frac{1}{N}\sum_{\alpha=1}^{N-1}\frac{1}{\mu_{\alpha}^{(1)}-z}\bigg{|}&=&\bigg{|}\int\frac{dF_{N}(x)}{x-z}
-(1-\frac{1}{N})\int\frac{dF_{N}^{(1)}(x)}{x-z}\bigg{|}\nonumber\\
&=&\frac{1}{N}\bigg{|}\int\frac{NF_{N}(x)-(N-1)F_{N}^{(1)}(x)}{(x-z)^{2}}dx\bigg{|}\nonumber\\
&\leq&\frac{1}{N}\int\frac{dx}{(x-z)^2}
=\frac{\pi}{N},\label{4.17}
\end{eqnarray}
where $F_{N}^{(1)}(x)$ is the empirical spectral distribution of $H^{(1)}$. In the inequality above, we use the well known interlacing property between eigenvalues of an Hermitian matrix and its submatrix as:
\begin{eqnarray*}
\lambda_{1}\leq\mu_{1}^{(1)}\leq\lambda_{2}\leq\cdots\leq\mu_{N-1}^{(1)}\leq\lambda_{N}.\label{4.18}
\end{eqnarray*}
Combining Lemma$~3.1$, (\ref{4.4.5}), (\ref{4.16.5}) and (\ref{4.17}) we conclude the proof of (i) of Lemma$~3.2$.

It is similar to get (\ref{3.4}) for $\widetilde{G}_{11}$, we now turn to
\begin{eqnarray*}
\widetilde{H}=YY^T-\sqrt{\frac M N}I_M.
\end{eqnarray*}
If we denote the $k$-th row of $Y_N$ by $\tilde{y}_k$ and the $(M-1)\times N$ matrix consisting of the other $M-1$ rows of $Y_N$ by $D_k$. And further introduce the matrix
\begin{eqnarray*}
\widetilde{H}^{(k)}=[\widetilde{H}^{(k)}_{jk}]_{N\times N}=:D_k^TD_k-\sqrt{\frac MN}I_N
\end{eqnarray*}
with its eigenvalues $\tilde{\mu}_{\alpha}^{(k)}, \alpha=1,\cdots,N$.
Then we have the following representation
$$
\widetilde{H}=\left(
\begin{array}{ccccc}
 \tilde{y}_{1}\cdot \tilde{y}_{1}-\sqrt{\frac{M}{N}} &\tilde{y}_1D_1^T  \\
 D_1\tilde{y}_1^T &D_{1}D_1^T-\sqrt{\frac{M}{N}}I_{M-1}  \\
\end{array}
\right).
$$
If we use the $\widetilde{E}_1\{\cdot\}$ to denote the expectation with respect to $\tilde{y}_1$, and set
\begin{eqnarray*}
\widetilde{V}=\tilde{y}_1D_1^TD_1\left(D_1^TD_1-\sqrt{\frac MN}-z\right)^{-1}\tilde{y}_1^T.\label{4.19}
\end{eqnarray*}
We have
\begin{eqnarray*}
\widetilde{G}_{11}&=&\frac{1}{\tilde{y}_1\cdot \tilde{y}_1-\sqrt{\frac MN}-z-\widetilde{V}}\nonumber\\
&=:&\frac{1}{-\sqrt{\frac MN}-z-m_N(z)+\tilde{r}_1},\label{4.20}
\end{eqnarray*}
where
\begin{eqnarray*}
\tilde{r}_1&=&\tilde{y}_1\cdot \tilde{y}_1+m_N(z)+\widetilde{V}\nonumber\\
&=&\tilde{y}_1\cdot \tilde{y}_1+m_N(z)-\frac1N\sum_{\alpha=1}^{N}\frac{\sqrt{\frac NM}\tilde{\mu}_{\alpha}^{(1)}+1}{\tilde{\mu}_{\alpha}^{(1)}+z}+\widetilde{V}-\widetilde{E}_1\{\widetilde{V}\}.\label{4.21}
\end{eqnarray*}
Now we denote the set of event
\begin{eqnarray*}
\Omega_\circ=:\bigg{\{}|f_N(z)-m_N(z)+\tilde{r}_1(z)|\geq\frac12\sqrt{\frac MN}\bigg{\}},
\end{eqnarray*}
and observe that on $\Omega_\circ^c$ there exists
\begin{eqnarray*}
\bigg{|}\bigg{(}\sqrt{\frac{M}{N}}+z+m_N(z)-\tilde{r}_1(z)\bigg{)}\bigg{(}\sqrt{\frac MN}+z+f_N(z)\bigg{)}\bigg{|}\geq C\frac MN .
\end{eqnarray*}
Then we have
\begin{eqnarray}
&&\mathbb{E}\bigg{|}\widetilde{G}_{11}+\frac{1}{\sqrt{\frac MN}+z+f_N(z)}\bigg{|}^2\nonumber\\
&=&\mathbb{E}\bigg{|}\frac{f_N(z)-m_N(z)+\tilde{r}_1(z)}{(\sqrt{\frac{M}{N}}+z+m_N(z)-\tilde{r}_1(z))(\sqrt{\frac MN}+z+f_N(z))}\bigg{|}^2\nonumber\\
&\leq&C\bigg{[}(\frac NM)^2+\frac{N}{M}\mathbb{P}(\Omega_\circ)\bigg{]}\cdot\mathbb{E}|f_N(z)-m_N(z)+\tilde{r}_1(z)|^2, \label{4.21.4}
\end{eqnarray}
where we used the fact that
\begin{eqnarray}
\bigg{|}\bigg{(}\sqrt{\frac{M}{N}}+z+m_N(z)-\tilde{r}_1(z)\bigg{)}\bigg{(}\sqrt{\frac MN}+z+f_N(z)\bigg{)}\bigg{|}\geq C\sqrt{\frac MN}\label{4.21.3}
\end{eqnarray}
holds on the full set $\Omega$. To see this, We denote the unit eigenvector of $\widetilde{H}^{(1)}$ corresponding to the eigenvalue $\tilde{\mu}_{\alpha}^{(1)}$ by $\tilde{v}_{\alpha}^{(1)}=(\tilde{v}_{\alpha}^{(1)}(1),\cdots,\tilde{v}_{\alpha}^{(1)}(N))$ $(\alpha=1,\cdots,N)$ and set
$$\tilde{\xi}_{\alpha}^{(1)}=\sqrt{MN}|\tilde{y}_{1}\cdot \tilde{v}_{\alpha}^{(1)}|^2.$$
Then we have
\begin{eqnarray}
\sqrt{\frac{M}{N}}+z+m_N(z)-\tilde{r}_1(z)&=&\sqrt{\frac{M}{N}}+z-\tilde{y}_1\bigg{[}I-(\tilde{H}^{(1)}+\sqrt{\frac MN})(\tilde{H}^{(1)}-z)^{-1}\bigg{]}\tilde{y}_1^T\nonumber\\
&=&\sqrt{\frac{M}{N}}+z+\frac{1}{\sqrt{MN}}\sum_{\alpha=1}^{N}\bigg{(}\frac{\sqrt{\frac MN}+\tilde{\mu}_{\alpha}^{(1)}}{\tilde{\mu}_{\alpha}^{(1)}-z}-1\bigg{)}\tilde{\xi}_{\alpha}^{(1)}\nonumber\\
&=&(\sqrt{\frac{M}{N}}+z)\left(1+\frac{1}{\sqrt{MN}}\sum_{\alpha=1}^{N}\frac{1}{\tilde{\mu}_{\alpha}^{(1)}-z}\tilde{\xi}_{\alpha}^{(1)}\right)\nonumber\\
&=:&(\sqrt{\frac{M}{N}}+z)(1+S).\label{4.21.6}
\end{eqnarray}
Taking into account that $\tilde{\mu}_{\alpha}^{(1)}\geq-\sqrt{\frac MN}$ and $\tilde{\xi}_{\alpha}^{(1)}\geq 0$, we have
\begin{eqnarray*}
\mathrm{Re}S\geq -C\sqrt{\frac MN}\mathrm{Im}S,\qquad \mathrm{Im}S\geq 0,
\end{eqnarray*}
which implies
\begin{eqnarray}
|1+S|\geq C\sqrt{\frac NM}.\label{4.21.7}
\end{eqnarray}
Then (\ref{4.21.3}) is a direct consequence of (\ref{4.21.6}) and (\ref{4.21.7}).

Now we proceed to the estimate of (\ref{4.21.4}). By observing
\begin{eqnarray*}
\mathbb{P}(\Omega_\circ)\leq \frac{4N}{M}\mathbb{E}|f_N(z)-m_N(z)+\tilde{r}_1(z)|^2,\nonumber
\end{eqnarray*}
it suffices to provide
\begin{eqnarray*}
\mathbb{E}|f_N(z)-m_N(z)+\tilde{r}_1(z)|^2\leq C_1\frac NM+C_2\frac1N\nonumber,
\end{eqnarray*}
which can be derived similarly as what we have done to (\ref{4.16}).
\end{proof}

\section{CLT For Linear Eigenvalue Statistics}
 To prove Theorem$~1.2$, we will follow the recent article \cite{SM} by Shcherbina. For there are only some technical differences, we will only state the main body of the proof in this section, and left all the technical details to the Appendix.
\begin{proof}[\emph{Proof of Theorem$~1.2$}] Similar to \cite{SM}, as we will see, $\omega_8$ is needed in our proof. So firstly we need to truncate the random variable $X_{jk}$ at $(MN)^{1/4}$, and then re-centralize it. To use the truncated matrix, it is necessary to show its linear eigenvalue statistics have the same limit distribution as the original one at first. If we denote the truncated matrix as
\begin{eqnarray*}
&&\widehat{Y}=[\widehat{Y}_{ij}]_{M\times N},\qquad \widehat{Y}_{ij}=Y_{ij}\mathbf{1}_{\{|Y_{ij}|\leq 1\}},\nonumber\\
&&\breve{Y}=[\breve{Y}_{ij}]_{M\times N}, \qquad  \breve{Y}_{ij}=\widehat{Y}_{ij}-\mathbb{E}\widehat{Y}_{ij}.\label{5.1}
\end{eqnarray*}
And further introduce
\begin{eqnarray*}
&&\widehat{H}_N=\widehat{Y}^T\widehat{Y}-\sqrt{\frac MN}I_N,\qquad \widehat{\mathcal{L}}_N[\varphi]=Tr\varphi(\widehat{H}),\nonumber\\
&&\breve{H}_N=\breve{Y}^T\breve{Y}-\sqrt{\frac MN}I_N,\qquad \breve{\mathcal{L}}_N[\varphi]=Tr\varphi(\breve{H}),\nonumber\\
&&\widehat{H}_N(t)=[\widehat{Y}+t(Y-\widehat{Y})]^T[\widehat{Y}+t(Y-\widehat{Y})]-\sqrt{\frac{M}{N}}I_{N}.\label{5.2}
\end{eqnarray*}
We denote the eigenvalues and corresponding eigenvectors of $\widehat{H}_N(t)$ by $\lambda_{i}(t)$ and $v_i(t)$ ($i=1,\cdots, N$) below. It suffices to prove that
\begin{eqnarray}
&&\mathbb{E}\{|e^{ix\mathcal{L}_N^\circ[\varphi]}-e^{ix\widehat{\mathcal{L}}_N^\circ[\varphi]}|\}\leq 2\mathbb{P}\{\widehat{Y}\neq Y\}+|x||\mathbb{E}\{\mathcal{L}_N[\varphi]\}-\mathbb{E}\{\widehat{\mathcal{L}}_N[\varphi]\}|\rightarrow 0,\label{5.3}\\
&&\mathbb{E}\{|e^{ix\widehat{\mathcal{L}}_N^\circ[\varphi]}-e^{ix\breve{\mathcal{L}}_N^\circ[\varphi]}|\}\leq2|x|
\mathbb{E}\{|\widehat{L}_N[\varphi]-\breve{\mathcal{L}}_N[\varphi]|\}\rightarrow 0.\label{5.4}
\end{eqnarray}
as $N$ tends to infinity. The proof of (\ref{5.4}) is similar to (\ref{5.3}), we only give the proof of (\ref{5.3}) below. $\mathbb{P}\{\widehat{Y}\neq Y\}$ is obviously $o(1)$ by the truncation. Also, we have
\begin{eqnarray}
&&\mathbb{E}\{|\mathcal{L}_N[\varphi]-\widehat{\mathcal{L}}_N[\varphi]|\}\nonumber\\
&&=\int_0^1dt\mathbb{E}\{\sum_i|\varphi'(\lambda_i(t))\lambda'_i(t)|\}\nonumber\\
&&\leq||\varphi'||_{\infty}\int_0^1dt\mathbb{E}\{\sum_{i}v_i^T(t)\widehat{H}'(t)v_i(t)\}\nonumber\\
&&\leq||\varphi'||_{\infty}\int_0^1dt\mathbb{E}\{Tr|(Y-\widehat{Y})^T\widehat{Y}+\widehat{Y}^T(Y-\widehat{Y})+2t(Y-\widehat{Y})^T(Y-\widehat{Y})|\}\nonumber\\
&&\leq||\varphi'||_{\infty}\int_0^1dt\mathbb{E}\{\sum_{i,j}|[(Y-\widehat{Y})^T\widehat{Y}+\widehat{Y}^T(Y-\widehat{Y})
+2t(Y-\widehat{Y})^T(Y-\widehat{Y})]_{ij}|\}\nonumber\\
&&\leq2||\varphi'||_{\infty}\int_0^1dt\sum_{i,j,k}(\mathbb{E}\{|(Y_{ki}-\widehat{Y}_{ki})|\cdot|\widehat{Y}_{kj}|\}+t\mathbb{E}\{|(Y_{ki}
-\widehat{Y}_{ki})(Y_{kj}-\widehat{Y}_{kj})|\})\nonumber\\
&&\leq2||\varphi'||_{\infty}\int_0^1dt\bigg{\{}\sum_{i,k}[\mathbb{E}\{|Y_{ki}-\widehat{Y}_{ki}|\cdot|\widehat{Y}_{ki}|\}+t\mathbb{E}\{|Y_{ki}-\widehat{Y}_{ki}|^2\}]\nonumber\\
&&+\sum_{i\neq j}\sum_{k}[\mathbb{E}|Y_{ki}-\widehat{Y}_{ki}|\cdot\mathbb{E}|\widehat{Y}_{kj}|
+t\mathbb{E}|Y_{ki}-\widehat{Y}_{ki}|\cdot\mathbb{E}|Y_{kj}-\widehat{Y}_{kj}|]\bigg{\}}.\label{5.5}
\end{eqnarray}
By the truncation, we have
\begin{eqnarray*}
&&\mathbb{E}|\widehat{Y}_{kj}|\leq C(MN)^{-\frac14},\nonumber\\
&&\mathbb{E}|Y_{ki}-\widehat{Y}_{ki}|\leq(MN)^{-\frac14}\int|X_{ki}|\mathbf{1}_{\{|X_{ki}|\geq(MN)^{\frac14}\}}dP\leq C(MN)^{-\frac54}.\nonumber\\
&&\mathbb{E}|Y_{ki}-\widehat{Y}_{ki}|^2\leq(MN)^{-\frac12}\int|X_{ki}|^2\mathbf{1}_{\{|X_{ki}|\geq(MN)^{\frac14}\}}dP\leq C(MN)^{-\frac54}.\label{5.6}
\end{eqnarray*}
Then it is easy to check (\ref{5.5}) tends to $0$ as $N$ tends to infinity, so (\ref{5.3}) holds.

So without loss of generality, we may and do assume in the sequel
\begin{eqnarray}
&&\mathbb{E}\{X_{jk}\}=0, \quad \mathbb{E}\{X_{jk}^2\}=1+o(1),\quad\mathbb{E}\{X_{jk}^4\}=\omega_4+o(1),\label{5.7}\\
&&\mathbb{E}\{X_{jk}^6\}\leq C(MN)^{\frac14}, \quad\mathbb{E}\{X_{jk}^8\}\leq C(MN)^{\frac34}.\label{5.8}
\end{eqnarray}
Relying on (\ref{5.7}) and (\ref{5.8}) we have the following lemma which collects all estimates we need to derive the central limit
theorem for  $\mathcal{L}_N[\varphi]$. \\\\
$\bf{Lemma~5.1}$ \emph{Under the assumptions (\ref{5.7}), (\ref{5.8}) and our basic assumption $N/M\rightarrow 0$ as $N\rightarrow\infty$, we have for any fixed $z$ with $\eta>0$ and $z_1, z_2: \mathrm{Im}z_{1,2}>0$, $1>\delta>0$ the following estimates hold:
\begin{eqnarray}
&&\mathrm{Var}\{\gamma_N\}\leq CN^{-1}\sum_{i=1}^{N}(\mathbb{E}\{|G_{ii}(z)|^{1+\delta}\})(\eta^{-3-\delta}+\eta^{-3+\delta}),\label{5.9}\\
&&\mathbb{E}\{(V^\circ)^3\},\quad\mathbb{E}\{(U^\circ)^3\},\quad \mathbb{E}\{|V^\circ|^4\},\quad\mathbb{E}\{|U^\circ|^4\}=O(N^{-\frac32}),\label{5.10}\\
&&N\mathbb{E}_1\{V^\circ(z_1)V^\circ(z_2)\}=\omega_4-1-\frac{\omega_4-1}{M}Tr[M^{[1]}(z_1)+M^{[1]}(z_2)]\nonumber\\
&&+\frac{\omega_4-3}{M}\sum_{i}M^{[1]}_{ii}(z_1)M^{[1]}_{ii}(z_2)+\frac{2}{M}Tr[M^{[1]}(z_1)M^{[1]}(z_2)]\nonumber\\
&&+\frac{1}{M}(\sqrt{\frac{M}{N}}+z_1)(\sqrt{\frac{M}{N}}+z_2)\gamma_N^{\circ(1)}(z_1)\gamma_N^{\circ(1)}(z_2),\label{5.11}\\
&&\mathrm{Var}\{N\mathbb{E}_1\{V^\circ(z_1)V^\circ(z_2)\}\},\quad \mathrm{Var}\{N\mathbb{E}_1\{V^\circ(z_1)U^\circ(z_2)\}\}=O\left(\frac NM\vee \frac1N\right),\label{5.12}\\
&&\mathbb{E}\{|\gamma_N^{\circ}|^4\}\leq CN^{1/2}\eta^{-12}+o(N^{1/2}),\quad\mathbb{E}\{|\gamma^{\circ(1)}_N-\gamma^\circ_N|^4\}=O(N^{-3/2}),\label{5.13}\\
&&|\mathbb{E}\{\gamma_{N}^{(1)}\}/N-f(z)|,\quad |\mathbb{E}^{-1}\{V\}+f(z)|=o(1).\label{5.14}
\end{eqnarray}}
We postpone the proof of Lemma$~5.1$ to the Appendix. Here we need to mention that the constants in (\ref{5.9}) or (\ref{5.13}) are independent with $\eta$.

By Proposition$~1$, Proposition$~3$ of \cite{SM} and (\ref{5.9}) we only need to prove Theorem$~2$ for the functions $\varphi=\varphi_{\eta}$ which are the convolution of some $\varphi_0$ with the Poisson kernel $P_{\eta}=\eta/\pi(x^2+\eta^2)$. And $\varphi_0$ is restricted to satisfy $\int|\varphi_0(\lambda)|d\lambda\leq C$ with some positive constant $C$. For such test function $\varphi$ we have
\begin{eqnarray*}
\mathcal{L}_N^\circ[\varphi]=\frac1\pi\int\varphi_0(\mu)\mathrm{Im}\gamma_N^\circ(z_\mu)d\mu,\quad z_\mu=\mu+i\eta.\label{5.15}
\end{eqnarray*}
Follow the notations of \cite{SM}, we set
\begin{eqnarray}
Z_N(x)=\mathbb{E}\{e^{ix\mathcal{L}_N^\circ[\varphi]}\},\quad e(x)=e^{ix\mathcal{L}_N^\circ[\varphi]}, \quad e_1(x)=e^{ix(\mathcal{L}^{(1)}_N[\varphi])^\circ},\label{5.16}
\end{eqnarray}
where $\mathcal{L}^{(1)}_N[\varphi]$ stands for the corresponding linear eigenvalue statistics of $H^{(1)}$.
Observe that Theorem$~2$ can be proved by providing that
\begin{eqnarray}
\frac{d}{dx}Z_N(x)=-x V[\varphi_0,\eta]Z_N(x)+o(1).\label{5.17}
\end{eqnarray}
To do this, we introduce
\begin{eqnarray}
Y_N(z,x)=:\mathbb{E}\{Tr G(z)e^\circ(x)\}=\sum_{k}\mathbb{E}\{G_{kk}(z)e^\circ(x)\}.\label{5.18}
\end{eqnarray}
Thus we have
\begin{eqnarray*}
\frac{d}{dx}Z_N(x)=\frac1{2\pi}\int\varphi_0(\mu)(Y(z_\mu,x)-Y(\bar{z}_\mu,x))d\mu.\label{5.19}
\end{eqnarray*}

We only need to deal with the first term in the summation (\ref{5.18}), the others are analogous,
so we can deal with $N\mathbb{E}\{V^{-1}e^\circ(x)\}:=T_1+T_2$ instead, where
\begin{eqnarray*}
&&T_1=-N\mathbb{E}\{(V^{-1})^\circ e_1(x)\},\label{5.20}\\
&&T_2=-N\mathbb{E}\{(V^{-1})^\circ(e(x)-e_1(x))\}.\label{5.21}
\end{eqnarray*}
For the calculations towards $T_1$ and $T_2$ are similar to those in the case studied in \cite{SM}, we will not state the tedious process. Indeed, by inserting the estimate in Lemma$~5.1$, we easily obtain
\begin{eqnarray*}
T_1=f^2(z)Y_N(z,x)+O\bigg{(}\frac{1}{\sqrt{N}}\vee\sqrt{\frac NM}\bigg{)},\label{5.22}
\end{eqnarray*}
and
\begin{eqnarray*}
T_2=ixZ_N(x)\int d\mu\varphi_0(\mu)\frac{D(z,z_\mu)-D(z,\bar{z}_\mu)}{2i\pi}+o(1),\label{5.23}
\end{eqnarray*}
where
\begin{eqnarray}
D(z,z_\mu)
&:=&2f^2(z)f^2(z_\mu)(1+f'(z_\mu))\bigg{(}\frac{f(z)-f(z_\mu)}{z-z_\mu}+\frac{\omega_4-1}{2}\bigg{)}\nonumber\\
&~&+2f^2(z)f(z_\mu)\frac{d}{dz_\mu}\bigg{(}\frac{f(z)-f(z_\mu)}{z-z_\mu}\bigg{)}.\label{5.24}
\end{eqnarray}
So we have
\begin{eqnarray*}
Y_N(z,x)=f^2(z)Y_N(z,x)+ixZ_N(x)\int d\mu\varphi_0(\mu)\frac{D(z,z_\mu)-D(z,\bar{z}_\mu)}{2i\pi(1-f^2(z))}+o(1).\label{5.25}
\end{eqnarray*}
Thus the variance in (\ref{5.17}) can be represented by
\begin{eqnarray*}
V[\varphi_0,\eta]&=&\frac{1}{4\pi^2}\int\int\varphi_0(\mu_1)\varphi_0(\mu_2)\bigg{(}C(z_{\mu_1},\bar{z}_{\mu_2})+C(\bar{z}_{\mu_1},z_{\mu_2})\nonumber\\
&&-C(z_{\mu_1},z_{\mu_2})-C(\bar{z}_{\mu_1},\bar{z}_{\mu_2})\bigg{)}d\mu_1d\mu_2,\label{5.26}
\end{eqnarray*}
where
\begin{eqnarray*}
C(z,z_\mu)=\frac{D(z,z_{\mu})}{1-f^2(z)}.\label{5.27}
\end{eqnarray*}
We can get (\ref{1.9}) by comparing (\ref{5.24}) to $(2.40)$ of \cite{SM}, then we complete the proof.
\end{proof}

\section{Appendix}
\begin{proof}[\emph{Proof of Lemma$~5.1$}]
We begin with (\ref{5.9}), we will use (\ref{4.1}) and (\ref{4.3}), and so we need to estimate (\ref{4.4}) more carefully. According to (\ref{4.9}) and (\ref{4.10}), we have
\begin{eqnarray}
\mathbb{E}_1\bigg{\{}\bigg{|}\frac{V-\mathbb{E}_1\{V\}}{\mathbb{E}_1\{V\}}\bigg{|}^2\bigg{\}}\leq C_1\frac{1}{MN}\frac{Tr|M^{[1]}|^2}{|z+\frac{1}{\sqrt{MN}}TrM^{[1]}|^2}+C_2\frac{1}{N}\frac{1}{|\mathbb{E}_1\{V\}|^2}\label{6.1}
\end{eqnarray}
and
\begin{eqnarray*}
\mathbb{E}_1\bigg{\{}\bigg{|}\frac{U-\mathbb{E}_1\{U\}}{\mathbb{E}_1\{V\}}\bigg{|}^2\bigg{\}}=C\frac{1}{MN}\frac{Tr|M^{[2]}|^2}{|z+\frac{1}{\sqrt{MN}}Tr M^{[1]}|^2}.\label{6.2}
\end{eqnarray*}
For (\ref{6.1}), if we set $N^{(1)}=B_1^TB_1$, we have
\begin{eqnarray*}
TrM^{[1]}=TrG^{(1)}N^{(1)}, \qquad Tr|M^{[1]}|^2=TrG^{(1)}N^{(1)}G^{(1)*}N^{(1)}.\label{6.3}
\end{eqnarray*}
So
\begin{eqnarray*}
&&\frac{1}{MN}Tr|M^{[1]}|^2\nonumber\\
&&=\frac{1}{MN}TrG^{(1)}N^{(1)}G^{(1)*}N^{(1)}\nonumber\\
&&\leq\frac{1}{MN}[TrG^{(1)}N^{(1)}G^{(1)*}]^{1-\delta}[TrG^{(1)}(N^{(1)})^{\frac{1+\delta}{\delta}}G^{(1)*}]^\delta\nonumber\\
&&\leq\frac{1}{MN}[TrG^{(1)}N^{(1)}G^{(1)*}]^{1-\delta}[Tr(N^{(1)})^{\frac{1+\delta}{\delta}}]^\delta\eta^{-2\delta}\nonumber\\
&&\leq\frac{1}{N}\bigg{[}\frac{1}{\sqrt{MN}}TrG^{(1)}N^{(1)}G^{(1)*}\bigg{]}^{1-\delta}\bigg{[}\frac{1}{N}\bigg{(}\frac{N}{M}\bigg{)}^{\frac{1+\delta}{2\delta}}
Tr(N^{(1)})^{\frac{1+\delta}{\delta}}\bigg{]}^\delta\eta^{-2\delta}\nonumber\\
&&\leq\frac{1}{N}\bigg{[}\frac{1}{\sqrt{MN}}TrG^{(1)}N^{(1)}G^{(1)*}\bigg{]}^{1-\delta}\bigg{[}\frac{1}{N}
Tr\bigg{(}\frac1MX^{(1)T}X^{(1)}\bigg{)}^{\frac{1+\delta}{\delta}}\bigg{]}^\delta\eta^{-2\delta}.\label{6.4}
\end{eqnarray*}
 Furthermore we have
\begin{eqnarray*}
\mathrm{Im}TrG^{(1)}N^{(1)}=\eta TrG^{(1)}N^{(1)}G^{(1)*},\label{6.5}
\end{eqnarray*}
which implies
\begin{eqnarray}
&&\mathbb{E}\bigg{\{}\frac{1}{MN}\frac{Tr|M^{[1]}|^2}{\big{|}z+\frac{1}{\sqrt{MN}}TrM^{[1]}\big{|}^2}\bigg{\}}\nonumber\\
&&\leq\frac{1}{N}\eta^{-1-\delta}\mathbb{E}\bigg{\{}\frac{\big{[}\frac{1}{N}
Tr\left(\frac1MX^{(1)T}X^{(1)}\right)^{\frac{1+\delta}{\delta}}\big{]}^\delta}{|\mathbb{E}_1\{V\}|^{1+\delta}}\bigg{\}}\nonumber\\
&&\leq\frac{1}{N}\eta^{-1-\delta}\mathbb{E}\bigg{\{}\big{[}\frac{1}{N}
Tr\left(\frac1MX^{(1)T}X^{(1)}\right)^{\frac{1+\delta}{\delta}}\big{]}^\delta \mathbb{E}_1|G_{11}|^{1+\delta}\bigg{\}}.\label{6.6}
\end{eqnarray}
Here we have used Jensen's inequality $|\mathbb{E}_1\{V\}|^{-1}\leq\mathbb{E}_1\{|V|^{-1}\}$.

If we control (\ref{6.6}) by using Theorem$~5.9$ of \cite{BaiS} by replacing $y$ with $N/M$, we can get
\begin{eqnarray*}
\mathbb{E}_1\left\{\left|\frac{V-\mathbb{E}_1\{V\}}{\mathbb{E}_1\{V\}}\right|^2\right\}\leq
C\frac1N \mathbb{E}|G_{11}|^{1+\delta}(\eta^{-1-\delta}+\eta^{-1+\delta}).
\end{eqnarray*}
Then by the fact that $\mathrm{Im}V=-\eta U-\eta$, we also have
\begin{eqnarray*}
\mathbb{E}\{|\frac{V-\mathbb{E}_1\{V\}}{\mathbb{E}_1\{V\}}\frac{U}{V}|^2\}\leq C\frac{1}{N}(\mathbb{E}\{|G_{11}|^{1+\delta}\})(\eta^{-3-\delta}+\eta^{-3+\delta})\label{6.7}
\end{eqnarray*}
and
\begin{eqnarray*}
\mathbb{E}\{|\frac{U-\mathbb{E}_1\{U\}}{\mathbb{E}_1\{V\}}|^2\}\leq C\frac1N(\mathbb{E}\{|G_{11}|^{1+\delta}\})(\eta^{-3-\delta}+\eta^{-3+\delta}).\label{6.8}
\end{eqnarray*}

Next we turn to the proof of the first inequality of (\ref{5.13}), which will be used in the proof of (\ref{5.10}). We use the following inequality for martingale (see \cite{DFJ})
\begin{eqnarray*}
\mathbb{E}\{|\gamma_N^\circ|^4\}\leq CN\sum_{k=1}^{N}\mathbb{E}\{|\gamma_N-\mathbb{E}_k\{\gamma_N\}|^4\}.\label{6.9}
\end{eqnarray*}
Similar to (\ref{4.3}), it suffices to check
\begin{eqnarray}
\mathbb{E}\{|U-\mathbb{E}_1\{U\}|^4\}\leq C N^{-\frac32}\eta^{-8}+o(N^{-\frac32})\label{6.10}
\end{eqnarray}
and
\begin{eqnarray}
\mathbb{E}\{|V-\mathbb{E}_1\{V\}|^4\}\leq C N^{-\frac32}\eta^{-4}+o(N^{-\frac32}).\label{6.11}
\end{eqnarray}
We only present the proof of (\ref{6.11}) below, (\ref{6.10}) is analogous.
Observing that $M^{[1]}_{ii}=(B_1G^{(1)}B_1^T)_{ii}$, (\ref{3.13}) and (\ref{3.15}) are valid as well. We have
\begin{eqnarray*}
\sum_{i}\mathbb{E}|M_{ii}^{[1]}|^4&\leq& C\eta^{-2}\frac MN\sum_{i}\mathbb{E}|M_{ii}^{[1]}|^2=CM\eta^{-4},\label{6.12}
\end{eqnarray*}
which implies
\begin{eqnarray*}
&&\mathbb{E}|V-\mathbb{E}_1\{V\}|^4\nonumber\\
&&\leq C\mathbb{E}\bigg{\{}\mathbb{E}_1\bigg{[}\bigg{|}\sum_{i\neq j}M^{[1]}_{ij}Y_{i1}Y_{j1}\bigg{|}^4+\bigg{|}\sum_{i}M_{ii}^{[1]}(Y_{i1}^2-\frac{1}{\sqrt{MN}})\bigg{|}^4+\bigg{|}y_1\cdot y_1-\sqrt{\frac{M}{N}}\bigg{|}^4\bigg{]}\bigg{\}}\nonumber\\
&&\leq C_1\frac{1}{(MN)^2}\mathbb{E}\{Tr|M^{[1]}|^4\}+C_2\frac{1}{(MN)^2}\sum_{i}\mathbb{E}\{|M_{ii}^{[1]}|^4\}\omega_8\nonumber\\
&&+C_3\frac{1}{(MN)^2}\mathbb{E}\{(Tr|M|^2)^2\}+C_4\frac{1}{MN^2}\omega_8\nonumber\\
&&\leq CN^{-\frac32}\eta^{-4}+o(N^{-\frac32}).\label{6.13}
\end{eqnarray*}

For (\ref{5.10}), we only deal with $\mathbb{E}\{|V^\circ|^4\}$ below, the others are similar.
\begin{eqnarray*}
V^\circ &=&V-\mathbb{E}_1\{V\}+\mathbb{E}_1\{V\}-\mathbb{E}\{V\}\nonumber\\
          &=&V-\mathbb{E}_1\{V\}+\frac{1}{\sqrt{MN}}(TrM^{[1]}-\mathbb{E}TrM^{[1]})\nonumber\\
          &=&V-\mathbb{E}_1\{V\}+\left(\frac{1}{N}+\frac{1}{\sqrt{MN}}z_1\right)(\gamma_{N}^{(1)}-\mathbb{E}\gamma_{N}^{(1)}).
\end{eqnarray*}
So use the first inequality of (\ref{5.13}) to $\gamma_N^{(1)}$, together with (\ref{6.13}) we can get (\ref{5.10}).

For (\ref{5.11}), note the following expansion:
\begin{eqnarray}
&&N\mathbb{E}_1\{V^\circ(z_1)V^\circ(z_2)\}\nonumber\\
&&=(\omega_4-1)-\frac{\omega_4-1}{M}Tr[M^{[1]}(z_1)+M^{[1]}(z_2)]\nonumber\\
&~&+\frac{\omega_4-3}{M}\sum_{i}M^{[1]}_{ii}(z_1)M^{[1]}_{ii}(z_2)+\frac{2}{M}Tr[M^{[1]}(z_1)M^{[1]}(z_2)]\nonumber\\
&~&+\frac{1}{M}[(TrM^{[1]}(z_1)-\mathbb{E}TrM^{[1]}(z_1))(TrM^{[1]}(z_2)-\mathbb{E}TrM^{[1]}(z_2))]\nonumber\\
&&=(\omega_4-1)-\frac{\omega_4-1}{M}Tr[M^{[1]}(z_1)+M^{[1]}(z_2)]\nonumber\\
&~&+\frac{\omega_4-3}{M}\sum_{i}M^{[1]}_{ii}(z_1)M^{[1]}_{ii}(z_2)+\frac{2}{M}Tr[M^{[1]}(z_1)M^{[1]}(z_2)]\nonumber\\
&~&+\frac{1}{M}(\sqrt{\frac{M}{N}}+z_1)(\sqrt{\frac{M}{N}}+z_2)\gamma_N^{\circ(1)}(z_1)\gamma_N^{\circ(1)}(z_2).\label{6.15}
\end{eqnarray}
To deal with the first estimate in (\ref{5.12}), we only need to take care of the variance of the third and fourth terms of (\ref{6.15}).
Using (\ref{3.15}) to $M^{[1]}_{ii}$ again we can get
\begin{eqnarray*}
\mathrm{Var}\{M^{[1]}_{ii}\}\leq C_1(\frac{N}{M})^2+C_2\frac1M,
\end{eqnarray*}
together with the trivial bound $|M^{[1]}_{ii}|\leq C\sqrt{\frac MN}$ we have
\begin{eqnarray*}
&&\mathrm{Var}\{\frac1M\sum_{i}M^{[1]}_{ii}(z_1)M^{[1]}_{ii}(z_2)\}\nonumber\\
&&\leq\frac1M\sum_{i}\mathrm{Var}\{M^{[1]}_{ii}(z_1)M^{[1]}_{ii}(z_2)\}\nonumber\\
&&\leq\frac 1M\sum_{i}\mathbb{E}\bigg{(}M^{[1]}_{ii}(z_1)M^{[1]\circ}_{ii}(z_2)+M^{[1]\circ}_{ii}(z_1)\mathbb{E}M^{[1]}_{ii}(z_2)\bigg{)}^2\nonumber\\
&&\leq C_1\frac NM+C_2\frac1N.
\end{eqnarray*}
Also we have
\begin{eqnarray}
&&TrM^{[1]}(z_1)M^{[1]}(z_2)\nonumber\\
&&=\sum_{\alpha=1}^{N-1}\frac{(\mu_{\alpha}^{(1)}+\sqrt{\frac{M}{N}})^2}{(\mu_{\alpha}^{(1)}-z_1)(\mu_{\alpha}^{(1)}-z_2)}\nonumber\\
&&=\bigg{[}\frac{M}{N}-z_1z_2+(\sqrt{\frac MN}+\frac{z_1+z_2}{2})(z_1+z_2)\bigg{]}\sum_{\alpha=1}^{N-1}\frac{1}{(\mu_{\alpha}^{(1)}-z_1)(\mu_{\alpha}^{(1)}-z_2)}\nonumber\\
&&+(\sqrt{\frac MN}+\frac{z_1+z_2}{2})\sum_{\alpha=1}^{N-1}\frac{2\mu_{\alpha}^{(1)}-z_1-z_2}{(\mu_{\alpha}^{(1)}-z_1)(\mu_{\alpha}^{(1)}-z_2)}+N-1\nonumber\\
&&=\bigg{[}\frac{M}{N}-z_1z_2+(\sqrt{\frac MN}+\frac{z_1+z_2}{2})(z_1+z_2)\bigg{]}Tr\frac{G^{(1)}(z_1)-G^{(1)}(z_2)}{z_1-z_2}\nonumber\\
&&+(\sqrt{\frac MN}+\frac{z_1+z_2}{2})Tr(G^{(1)}(z_1)+G^{(1)}(z_2))+N-1.\nonumber
\end{eqnarray}
Using (\ref{5.9}) to $TrG^{(1)}=\gamma_N^{(1)}$ we have
\begin{eqnarray*}
\mathrm{Var}\left\{\frac1M TrM^{[1]}(z_1)M^{[1]}(z_2)\right\}\leq \frac{1}{N^2}.
\end{eqnarray*}
 The estimate towards $\mathrm{Var}\{N\mathbb{E}_1\{V^\circ(z_1)U^\circ(z_2)\}\})$ is straightforward by taking into account $(G^{(1)})^2=dG^{(1)}/dz$.

 The second estimate of (\ref{5.13}) is a direct consequence of (\ref{4.2}) and (\ref{5.10}). The first part of (\ref{5.14}) is just the consequence of (\ref{1.0}) if we replace $f_N(z)=\gamma_N/N$ by $\gamma_N^{(1)}/N$, and second one follows directly from the fact
 $$\mathbb{E}\{V\}=z+\frac1N\gamma_N^{(1)}+O(\sqrt{\frac NM})=z+f(z)+o(1)$$
 and
 $$f(z)=-\frac{1}{z+f(z)}.$$
So we complete the proof.
\end{proof}


\begin{thebibliography}{00}

\bibitem{APS}
S. Albeverio, L. Pastur, M. Shcherbina.: \emph{On the 1/n expansion for some unitary invariant ensemble of random matrices.} Commun. Math. Phys. 224, 271-305 (2001).
\bibitem{AZ}
G.W. Anderson, O. Zeitouni: \emph{CLT for a band matrix model}. Probab. Theory and Related Fields. V.134, 283-338 (2006).
\bibitem{BaS}
Z.D. Bai, J.W. Silverstein: \emph{CLT for linear spectral statistics of large dimensional sample covariance matrices}. Ann. of Prob. V.32, 553-605 (2004).
\bibitem{BaiS}
Z.D. Bai, J.W. Silverstein: \emph{Spectral analysis of large dimensional random matrices} Mathematics Monograph Series 2, Science Press, Beijing.
\bibitem{BY}
Z.D. Bai, Y.Q. Yin: \emph{Convergence to the semicircle law}. Ann. of Prob. Volumn 16, No.2, 863-875 (1988).
\bibitem{DFJ}
S.W. Dharmadhikari, V. Fabian, K. Jogdeo: \emph{Bounds on the moments of martingales.} Ann. Math. Statist. v.39, p. 1719-1723 (1968).
\bibitem{EK}
N. EL Karoui: \emph{On the largest eigenvalue of Wishart matrices with identity covariance when $n, p$ and $p/n\rightarrow\infty$}. Preprint, arXiv: math.ST/0309355.
\bibitem{Hsu}
P.L. Hsu: \emph{On the distribution of roots of certain determinantal equations.} Ann. Eugenics 9, 250-258.
\bibitem{Joh}
K. Johansson: \emph{On fluctuations of eigenvalues of random Hermitian matrices}. Duke Math. J. V.91, 151-204 (1998).
\bibitem{KKP}
A.M. Khorunzhy, B.A. Khoruzhenko, L.A. Pastur: \emph{On asymptotic properties of large random matrices with independent entries.}
J. Math. Phys. 37, 5033 (1996).
\bibitem{KK}
A.M. Khorunzhy, W. Kirsch: \emph{On asymptotic expansions and scales of spectral universality in band random matrix ensembles.} Commun. Math. Phys 231, 223-255 (2002).
\bibitem{LP}
A. Lytova, L. Pastur: \emph{central limit theorem for linear eigenvalue statistics of random matrices with independent entries} Ann. of prob. 37 No.5, 1778-1840 (2009).
\bibitem{MP}
V.A. Mar\u{c}enko, L.A. Pastur: \emph{Distribution of eigenvalues for some sets of random matrices.} Math. USSR-Sb. 1, 457-483 (1967).
\bibitem{SM}
M. Shcherbina: \emph{Central limit theorem for linear eigenvalue statistics of the Wigner and sample covariance random matrices} Preprint, arXiv:1101.3249v1.
\bibitem{SS}
Ya. Sinai, A. Soshnikov: \emph{Central limit theorem for traces of large random symmetric matrices with independent matrix elements}. Bol.Soc.Brasil.Mat.(N.S.) V.29, 1-24 (1998).
\bibitem{WE}
E. Wigner: \emph{Characteristic vectors of bordered matrices with infinite dimensions.} Ann. of Math. 62 , 548-564 (1955).
\bibitem{WJ}
J. Wishart: \emph{The generalized product moment distribution in samples from a normal multivariate population.} Biometrika A 20, 32-43 (1928).

\end{thebibliography}
\end{document}